\documentclass{amsart}
\usepackage{amssymb}

\newtheorem{thm}{Theorem}[section]
\newtheorem{cor}[thm]{Corollary}
\newtheorem{conj}[thm]{Conjecture}
\newtheorem{corconj}[thm]{Conjecture-Corollary}
\newtheorem{prop}[thm]{Proposition}
\newtheorem{lem}[thm]{Lemma}

\newcommand\F{{\mathbb F}_2}
\newcommand\divs{{\big|\,}}
\newcommand\lcm{\operatorname{lcm}}
\newcommand\Gal{\operatorname{Gal}}
\newcommand\defeq{\overset{\text{\tiny def}}{=}}

\begin{document}
\bibliographystyle{amsplain}
\nocite{menezes-:ffields}
\title{Lamps, Factorizations and Finite Fields}
\author{Laurent Bartholdi}
\address{\parbox[b]{.4\linewidth}{Section de Math\'ematiques\\
    Universit\'e de Gen\`eve\\ CP 240, 1211 Gen\`eve 24\\
    Switzerland}}
\email{Laurent.Bartholdi@math.unige.ch}
\date{\today}
\maketitle

\section{Introduction}
The origin of this study is the 1993 International Mathematical
Olympiads, held at Istanbul, Turkey. Problem \# 6, which occurred on
day 2, reads:

\begin{quote}
  Let $n>1$ be an integer. There are $n$ lamps $L_0,\dots,L_{n-1}$
  arranged in a circle. Each lamp is either ON or OFF. A sequence of
  steps $S_0,S_1,\dots,S_j,\dots$ is carried out. Step $S_j$ affects
  the state of $L_j$ only (leaving the state of all others lamps
  unalterated) as follows:
  \begin{quote}
    if $L_{j-1}$ is ON, $S_j$ changes the state of $L_j$ from ON to
    OFF or from OFF to ON;\\
    if $L_{j-1}$ is OFF, $S_j$ leaves the state of $L_j$ unchanged. 
  \end{quote}
  The lamps are labeled mod $n$, that is,
  $$L_{-1}=L_{n-1},\quad L_0=L_n,\quad L_1=L_{n+1},\text{ etc.}$$
  Initially all lamps are ON. Show that
  \renewcommand\descriptionlabel[1]{\hspace\labelsep\textbf{(#1)}}
  \begin{description}
  \item[a] there is a positive integer $M(n)$ such that after $M(n)$
    steps all the lamps are ON again;
  \item[b] if $n$ has the form $2^k$ then all the lamps are ON after
    $n^2-1$ steps;
  \item[c] if $n$ has the form $2^k+1$ then all the lamps are ON after
    $n^2-n+1$ steps.
  \end{description}
\end{quote}

In this note we answer the Olympiads question using elementary algebra
over finite fields, and exhibit an interesting phenomenon when $n$ is
one less than a power of two. More generally, we are interested in the
minimal time $t(n)\ge1$ such that after repeating $t(n)$ times the
above instructions all lamps are again lit.

It turns out this question is tightly related to the factorization of
the polynomial $\Phi_n=X^n+X+1$ over the field $\F$. For $n=2^k$ or
$2^k+1$ it has only small factors, and there is a surprising
connection between the factorization of $\Phi_n$ and that of
$\Phi_{2^n-1}$.

Only undergraduate abstract algebra knowledge is assumed from the
reader; however unsolved problems appear, for instance in
Conjecture~\ref{conj:2^n-1}. It would be interesting to know by what
means pre-university students solved this Olympiad problem.

\section{An Algebraic Reformulation}
We let a lamp's state be represented by $0,1\in\F$ for unlit and lit
respectively, and number the lamps counterclockwise from $0$ to $n-1$
in such a way that we are about to alter the lamp at position $n-2$.
We denote by $(a_0,\dots,a_{n-1})$ the lamps' state. One step of
evolution amounts then to the following: replace $a_{n-2}$ by
$a_{n-2}+a_{n-1}$, and move to position $n-3$. The process is
invariant under rotation of the circle, so we may renumber the lamps so
that we are again at position $n-2$, and describe one step of
evolution as the operation
\begin{equation}\label{eq:evolution}
  (a_0,\dots,a_{n-1})\mapsto(a_{n-1},a_0,\dots,a_{n-3},a_{n-2}+a_{n-1}).
\end{equation}

In turn, the lamps' state $(a_0,\dots,a_{n-1})$ is conveniently encoded
as a polynomial,
\begin{equation}\label{eq:f}
  f = \sum_{i=0}^{n-1}a_iX^i \quad\in \F[X]/(X^n+X^{n-1}+1).
\end{equation}
The reason $f$ is represented as a polynomial in this peculiar ring is
that one step of evolution described in~(\ref{eq:evolution})
translates, in terms of polynomials, to the operation `$f:=X\cdot f$'.
Indeed the direct translation of~\ref{eq:evolution} in the polynomial
ring is
$$f\mapsto Xf + a_{n-1}(1+X^{n-1}+X^n);$$
or, in other words: the conversion from the list representation to the
polynomial one is linear; a lit lamp at position $i$ corresponds to
$X^i$, which evolves to $X^{i+1}$, and $i+1$ is the new position of
the lamp; and a lit lamp at position $n-1$ corresponds to $X^{n-1}$,
which evolves in $X^n=X^{n-1}+1$, which maps back to the lamp at
position $0$ and a switched lamp at position $n-1$.

Note now that $X$ is invertible: $1/X = X^{n-1}+X^{n-2}$. The ring
in~(\ref{eq:f}) is naturally isomorphic, via the map $X\mapsto X^{-1}$,
to
$$R_n\defeq\F[X]/(X^n+X+1).$$
We shall consider the evolution `$f:=X\cdot f$' occurring in
$R_n$; this amounts to consider the original question with time moving
backwards.

We denote by $R_n^\times$ the group of invertible elements of
$R_n$. The initial position corresponds to
$$f_0 = \sum_{i=0}^{n-1}X^i = \frac{X^n+1}{X+1} = X^{1-n} \in R_n^\times.$$
Thus $$t(n) = \min\{t\ge 1\,|\,X^t=1 \text{ in } R_n\} = |\langle X\rangle|,$$
where $\langle X\rangle$ is viewed as a subgroup (not an ideal!) of
$R_n^\times$. We have proved the
\begin{prop} $$t(n)<\infty.$$
  More precisely, $t(n)<2^n$, and divides $|R_n^\times|$.
\end{prop}
\begin{proof}
  $t(n)$ is the order of a subgroup of $R_n^\times$, and $R_n^\times$
  is a finite group of order at most $2^n-1$.
\end{proof}
We shall later give more details about the structure of $R_n^\times$;
for now explicit values of $t(n)$ can be given in a few special cases:
\begin{prop} If $n$ is a power of two (say $n=2^k$), then
  $$t(n) = n^2 - 1.$$
\end{prop}
\begin{proof} We compute
  $$X^{n^2} = (X^n)^n = (X+1)^n = X^n + 1 = X,$$
  so $X^{n^2-1} = 1$. Conversely, if $n\le t<n^2-n$, write $t=ni+j$ with
  $1\le i<n-1$ and $0\le j<n$, and note that the polynomial
  $X^t=X^{ni+j}=(X+1)^iX^j$ has degree at most $2n-3$ and span $i$;
  write it as $f+X^ng$ with $f$ and $g$ of degree less than $n$. It is
  equal, in $R_n$, to $f+(X+1)g$ where the two summands don't overlap,
  and therefore cannot equal $X$. If $2\le t<n$, it is clear that
  $X^t\neq X$, and if $n^2-n\le t<n^2$ the same holds by symmetry.
\end{proof}
\begin{prop} If $n$ is one more than a power of two (say $n=2^k+1$), then
  $$t(n) = n^2 - n + 1.$$
\end{prop}
\begin{proof} We compute
  $$X^{n^2} = (X^n)^n = (X+1)^n = X^n + X^{n-1} + X + 1 = X^{n-1},$$
  so $X^{n^2-n+1} = 1$. The argument in the proof of the previous
  proposition shows that no smaller value satisfies this equation.
\end{proof}

In case $k$ is one less than a power of two, say $k=2^n-1$, there is a
peculiar phenomenon:
\begin{prop}\label{prop:2^n-1}
  For all $n\ge2$,
  $$t(2^n-1)\divs2^{t(n)}-1.$$
\end{prop}
\begin{proof}
  In $R_{2^n-1}$, we may consider a subset
  $$Q_n = \{\sum_{i=0}^{n-1} a_i X^{2^i}\,|\,a_i\in\F\}.$$
  It is a
  vector subspace of dimension $n$, as the $X^i$ are linearly
  independent for $0\le i<2^n-1$. Elements of $R_{2^n-1}$ are
  polynomials and therefore can be composed, an operation we denote by
  $\circ$. This operation is internal to $Q_n$, and endows $Q_n$ with an
  $\F$-algebra structure: $f(g(X)+h(X))=f(g(X))+f(h(X))$ as soon as
  all the monomials of $f$ have degree a power of $2$. Moreover, $Q_n$
  is Abelian (on the basis $\{X^i\}$ we have $X^i\circ
  X^j=X^{ij}=X^j\circ X^i$), and $Q_n\cong R_n$ through the natural
  map $X^i\mapsto X^{2^i}$ extended by linearity.  Indeed
  $$X^{i+j} = X^i\cdot X^j\mapsto X^{2^i}\circ X^{2^j} = X^{2^{i+j}},$$
  and for any $f\in\F[X]$
  $$0 = f\cdot(X^n+X+1)\mapsto \widehat f\circ(X^{2^n} + X^2 + X),$$
  where $\widehat f$ is a polynomial divisible by $X$. It follows that
  any polynomial representing $0$ in $Q_n$ maps to a multiple (for
  $\cdot$) of $X^{2^n} + X^2 + X$, which in turn represents $0$ in
  $R_{2^n-1}$.

  Now the evolution `$f:= X\cdot f$' is mapped in $Q_n$ to
  `$g:= X^2\circ g = g^2$'; thus for all $t$ such that $X^t=1$
  in $R_n$, one has $X^{\circ t} = X$ in $Q_n$, and $X^{2^t-1} = 1$ in
  $R_{2^n-1}$.
\end{proof}

The following conjecture relies on numerical evidence. It has been
checked for $n\le16$ using \textsc{Gap}~\cite{gap:manual}
and~\textsc{Pari-GP}~\cite{pari:user} and their finite field algorithms.
\begin{conj}\label{conj:2^n-1}
  For all $n\ge2$,
  $$t(2^n-1)=2^{t(n)}-1.$$
\end{conj}

Recall that a polynomial $f\in\F[X]$ is \emph{primitive} if
$(\F[X]/f)^\times$ is generated by $X$.  A striking consequence of
Conjecture~\ref{conj:2^n-1} is the following
\begin{corconj}
  Let $n_0=2$ and define recursively $n_{i+1}=2^{n_i}-1$ for $i\ge0$.
  Then $X^{n_i}+X+1$ is irreducible and primitive in $\F[X]$ for all
  $i\ge0$.
\end{corconj}
\begin{proof}
  $\Phi_n$ is irreducible and primitive if and only if $R_n$ is a
  field and $R_n^\times$ is generated by $X$; this is equivalent to
  $t(n)=2^n-1$, its maximal possible value. We have
  $$t(n_{i+1})=2^{t(n_i)}-1=2^{2^{n_i}-1}-1=2^{n_{i+1}}-1,$$
  the first and third equalities following from
  Conjecture~\ref{conj:2^n-1} and the second from induction.
\end{proof}

\section{More Results on the Factorization of $X^n+X+1$}
We now turn to a more thorough study of the polynomial $\Phi_n =
X^n+X+1$ over $\F$. The behaviour of $t(n)$ is closely related to the
structure of the algebra $R_n$, which in turn is determined by the
factorization of $\Phi_n$.

We denote by $\chi$ the Frobenius
automorphism~\cite[page~9]{serre:arith} of $R_n$. Recall that any
algebra $A$ over $\F$ has an endomorphism defined by $\chi(g)=g^2$; if
$A$ is a finite field of degree $d$, then $\chi$ is invertible, of
order $d$, and generates the Galois group $\Gal(A/\F)$.

We show first that $\chi$ is invertible in $R_n$. For this purpose,
suppose $g\in\F[X]$ satisfies $g^2\equiv 0\mod\Phi_n$. It then
follows that $g\equiv 0\mod\Phi_n$, by the
\begin{lem}\label{lem:nrf}
  The $\Phi_n$ do not have repeated factors.
\end{lem}
\begin{proof}
  It suffices to show that $(\Phi_n,\Phi_n')=1$; if $n$ is even, then
  $\Phi_n'=1$, while if $n$ is odd, then $\Phi_n'=X^{n-1}+1$ and
  $$(\Phi_n,\Phi_n')\divs\Phi_n-X\Phi_n'=1.$$
\end{proof}

As a consequence, $R_n$ is semisimple, i.e.\ decomposes as a direct
sum of fields. Let $\Phi_n$ factor as $f_{n,1}\cdots f_{n,r_n}$, with
$f_{n,i}$ irreducible polynomials of degree $d_{n,i}$. Then $R_n$ splits
as
$$R_n = F_{n,1}\oplus\dots\oplus F_{n,r_n},$$
where the $F_{n,i}$ are field extensions of $\F$ of degree $d_{n,i}$. Note
then that the order of the Frobenius automorphism $\chi$ is $d_{n,i}$ in
$F_{n,i}$, and therefore is $\lcm\{d_{n,i}\}_{1\le i\le r_n}$ in $R_n$.
The following lemma is straightforward:
\begin{lem}\label{lem:lcm}
  For $i\in\{1,\dots,r_n\}$, let $\pi_i$ be the natural map
  $R_n\twoheadrightarrow F_{n,i}$. Then
  $$t(n) = \lcm\left\{|\langle\pi_i(X)\rangle|\right\}_{1\le i\le r}.$$
  In particular, $t(n)$ divides $u(n)=\lcm\{2^{d_{n,i}}-1\}_{1\le i\le
    r_n}$ (see Table~\ref{tbl:1}).
\end{lem}
\begin{table}
\caption{Factorizations of $X^n+X+1$ in $\F$ and corresponding $t(n)$
  and $u(n)$ (see Lemma~\ref{lem:lcm}}
\label{tbl:1}
\begingroup
\tiny
\begin{tabular}{|r|c|c|c|}
\hline
$n$     & $t(n)=|\langle X\rangle|$ & $u(n)/t(n)$ & $X^n + X + 1\quad(\hbox{\bf mod 2})$ \\ \hline
2       & 3     & 1     & $X^2 + X + 1$ \\
3       & 7     & 1     & $X^3 + X + 1$ \\
4       & 15    & 1     & $X^4 + X + 1$ \\
5       & 21    & 1     & $(X^2 + X + 1)(X^3 + X^2 + 1)$ \\
6       & 63    & 1     & $X^6 + X + 1$ \\
7       & 127   & 1     & $X^7 + X + 1$ \\
8       & 63    & 1     & $(X^2 + X + 1)(X^6 + X^5 + X^3 + X^2 + 1)$ \\
9       & 73    & 7     & $X^9 + X + 1$ \\
10      & 889   & 1     & $(X^3 + X + 1)(X^7 + X^5 + X^4 + X^3 + 1)$ \\
\hline
11      & 1533  & 1     & $(X^2 + X + 1)(X^9 + X^8 + X^6 + X^5 + X^3 + X^2 + 1)$ \\
12      & 3255  & 1     & $(X^3 + X^2 + 1)(X^4 + X^3 + 1)(X^5 + X^3 + X^2 + X + 1)$ \\
13      & 7905  & 1     & $(X^5 + X^4 + X^3 + X + 1)(X^8 + X^7 + X^5 + X^3 + 1)$ \\
14      & 11811 & 1     & $(X^2 + X + 1)(X^5 + X^3 + 1)(X^7 + X^6 + X^5 + X^2 + 1)$ \\
15      & 32767 & 1     & $X^{15} + X + 1$ \\
16      & 255   & 1     & $(X^8 + X^6 + X^5 + X^3 + 1)(X^8 + X^6 + X^5 + X^4 + X^3 + X + 1)$ \\
17      & 273   & 15    & $(X^2 + X + 1)(X^3 + X + 1)(X^{12} + X^{11} + X^{10} + X^9 + X^8 + X^6 + X^4 + X + 1)$ \\
18      & 253921 & 1    & $(X^5 + X^2 + 1)(X^{13} + X^{10} + X^8 + X^7 + X^4 + X^3 + X^2 + X + 1)$ \\
19      & 413385 & 1    & $(X^3 + X^2 + 1)(X^4 + X + 1)(X^5 + X^4 + X^2 + X + 1)(X^7 + X^5 + X^4 + X^3 + X^2 + X + 1)$ \\
20      & 761763 & 1    & $(X^2 + X + 1)(X^5 + X^4 + X^3 + X^2 + 1)(X^{13} + X^{11} + X^{10} + X^9 + X^7 + X^4 + 1)$ \\
\hline
21      & 5461  & 3     & $(X^7 + X^5 + X^3 + X + 1)(X^{14} + X^{12} + X^7 + X^6 + X^4 + X^3 + 1)$ \\
22      & 4194303 & 1   & $X^{22} + X + 1$ \\
23      & 2088705 & 1   & $(X^2 + X + 1)(X^8 + X^6 + X^3 + X^2 + 1)(X^{13} + X^{12} + X^{11} + X^8 + X^7 + X^5 + 1)$ \\
24      & 2097151 & 1   & $(X^3 + X + 1)(X^{21} + X^{19} + \dots + X^3 + 1)$ \\
25      & 10961685 & 1  & $(X^6 + X^5 + X^2 + X + 1)(X^8 + X^4 + X^3 + X^2 + 1)(X^{11} + X^{10} + X^9 + X^8 + X^7 + X^4 + 1)$ \\
26      & 298935 & 1    & $(X^2 + X + 1)(X^3 + X^2 + 1)(X^9 + X^7 + X^5 + X^4 + X^3 + X^2 + 1)(X^{12} + X^{10} + \dots + X^2 + 1)$ \\
27      & 125829105 & 1 & $(X^4 + X^3 + 1)(X^{23} + X^{22} + \dots + X + 1)$ \\
28      & 17895697 & 15 & $X^{28} + X + 1$ \\
29      & 402653181 & 1 & $(X^2 + X + 1)(X^{27} + X^{26} + \dots + X^2 + 1)$ \\
30      & 10845877 & 99 & $X^{30} + X + 1$ \\
\hline
31      & 2097151 & 1   & $(X^3 + X + 1)(X^7 + X^3 + 1)(X^{21} + X^{19} + X^{18} + X^{15} + X^{14} + X^{11} + X^8 + X^7 + X^5 + X^4 + 1)$ \\
32      & 1023  & 1     & $(X^2 + X + 1)(X^{10} + X^9 + X^8 + X^3 + X^2 + X + 1)$ \\
        & & & \hfill $(X^{10} + X^9 + X^8 + X^6 + X^5 + X + 1)(X^{10} + X^9 + X^8 + X^4 + X^3 + X^2 + 1)$ \\
33      & 1057  & 31 & $(X^3 + X^2 + 1)(X^{15} + X^{10} + X^9 + X^8 + X^4 + X^3 + X^2 + X + 1)$ \\
        & & & \hfill $(X^{15} + X^{14} + X^{13} + X^{11} + X^{10} + X^7 + X^6 + X^3 + 1)$ \\
34      & 255652815 & 21 & $(X^4 + X + 1)(X^{30} + X^{27} + \dots + X^4 + 1)$ \\
35      & 3681400539 & 7 & $(X^2 + X + 1)(X^{33} + X^{32} + \dots + X^2 + 1)$ \\
36      & 22839252821 & 3 & $(X^9 + X^7 + X^5 + X + 1)(X^{10} + X^7 + X^5 + X^3 + X^2 + X + 1)(X^{17} + X^{15} + \dots + X + 1)$ \\
37      & 137438167041 & 1 & $(X^{18} + X^{17} + \dots + X^5 + 1)(X^{19} + X^{18} + X^{17} + X^{13} + X^{12} + X^{10} + X^5 + X + 1)$ \\
38      & 25769803773 & 1 & $(X^2 + X + 1)(X^3 + X + 1)(X^{33} + X^{32} + \dots + X + 1)$ \\
39      & 178979337621 & 1 & $(X^6 + X^5 + X^4 + X + 1)(X^7 + X^6 + X^3 + X + 1)(X^{26} + X^{23} + \dots + X + 1)$ \\
40      & 320319056211 & 1 & $(X^3 + X^2 + 1)(X^{10} + X^8 + X^4 + X^3 + X^2 + X + 1)$ \\
        & & & $(X^{13} + X^9 + X^8 + X^6 + X^5 + X^4 + 1)(X^{14} + X^{13} + \dots + X^3 + 1)$ \\
\hline
41      & 545460846465 & 1 & $(X^2 + X + 1)(X^7 + X^6 + X^5 + X^4 + 1)(X^{32} + X^{30} + \dots + X^2 + 1)$ \\
42      & 1374389534715 & 1 & $(X^4 + X^3 + 1)(X^{38} + X^{37} + \dots + X + 1)$ \\
43      & 8521215115233 & 1 & $(X^5 + X^3 + X^2 + X + 1)(X^{38} + X^{36} + \dots + X^2 + 1)$ \\
44      & 12781822672803 & 1 & $(X^2 + X + 1)(X^5 + X^4 + X^3 + X + 1)(X^{37} + X^{35} + \dots + X + 1)$ \\
45      & 137434726401 & 1 & $(X^3 + X + 1)(X^5 + X^3 + 1)(X^{15} + X^{12} + X^4 + X^3 + 1)(X^{22} + X^{18} + \dots + X^3 + 1)$ \\
46      & 23456248059221 & 3 & $X^{46} + X + 1$ \\
47      & 1466015503701 & 3 & $(X^2 + X + 1)(X^3 + X^2 + 1)(X^{42} + X^{38} + \dots + X^5 + 1)$ \\
48      & 40209483820471 & 1 & $(X^{15} + X^{14} + X^{12} + X^{11} + X^{10} + X^6 + 1)(X^{33} + X^{32} + \dots + X + 1)$ \\
49      & 64677154575 & 17 & $(X^4 + X + 1)(X^5 + X^2 + 1)(X^{40} + X^{36} + \dots + X^2 + 1)$ \\
50      & 272662240182303 & 1 & $(X^2 + X + 1)(X^5 + X^4 + X^2 + X + 1)(X^{14} + X^{10} + \dots + X + 1)$ \\
        & & & $(X^{29} + X^{25} + X^{15} + X^{13} + X^8 + X^7 + 1)$ \\
\hline
\end{tabular}
\endgroup
\end{table}

\begin{prop}
  $\F[X]/\Phi_{2^k}\defeq A$ splits as a direct sum of fields of
  degree dividing $2k$.

  $\F[X]/\Phi_{2^k+1}\defeq B$ splits in factors of degree dividing $3k$.
\end{prop}
\begin{proof}
  Let us denote by $\psi:g\mapsto g^{2^k}$ the $k$-th power of the Frobenius
  automorphism. We must show that $\psi^2=1$ in $A$, and $\psi^3=1$ in
  $B$; but in $A$ we have $\psi(X)=1+X$ of order $2$, and in $B$ we
  have $\psi(X) = \frac{1+X}{X}$ of order $3$.
\end{proof}

This is in accordance with the results in the previous section:
$t(n=2^k)=n^2-1 = 2^{2k}-1$, and $t(n=2^k+1)=n^2-n+1\divs2^{3k}-1$.
Remark that the two transformations $X\mapsto 1+X$ and
$X\mapsto\frac{1+X}{X}$ of $R_n$ lift to
$PGL_2(\F)=\operatorname{Aut}(\F(X))$. These are the only possible
``systematic lifts'', and explains the special behaviour of $R_{2^k}$
and $R_{2^k+1}$.

For any polynomial $f=\sum a_i X^i\in\F[X]$, let us denote by
$\widehat f=\sum a_i X^{2^i}\in\F[X]$ the hat-polynomial of $f$.
(Sometimes $\widehat f$ is called a \emph{linearized polynomial} or a
\emph{$2$-polynomial}; see~\cite[\S3.4]{lidl-n:ffields}.)
Hat-polynomials can be multiplied, but also composed as in the proof
of Proposition~\ref{prop:2^n-1}. The composition operation $\circ$ is
linear thanks to the fact that all monomials in hat-polynomials have
degree a power of the field's characteristic, $2$; indeed
$$\widehat f\circ\widehat g=f(\chi)(\widehat g)=\widehat{f\cdot g}.$$
Let us note
$S_n\defeq\F[X]/\widehat{\Phi_n}=\F[X]/(X\cdot\Phi_{2^n-1})$; then for
any $f\in R_n$ we may naturally see $\widehat f\in
\widehat{R_n}\subset S_n$, and there is a natural embedding of $R_n$
in $\operatorname{End}(\widehat{R_n})$ given by $f\mapsto f(\chi)$,
with $f(\chi)(X)=\widehat f$. Note that under this embedding $X$ maps
to the Frobenius automorphism of $S_n$.

While $R_n$ decomposes as a direct sum, $S_n$ decomposes naturally as
a \emph{tensor product}. Recall that the tensor product of two
algebras $A$ and $B$ with bases $\{a_i\}$ and $\{b_j\}$
respectively is the algebra with basis $\{a_i\otimes b_j\}$ and
multiplication $(a_i\otimes b_j)(a_{i'}\otimes
b_{j'})=a_ia_{i'}\otimes b_jb_{j'}$. If $A=\F[X]/f(X)$ and
$B=\F[Y]/g(Y)$, one may take as bases $\{a_i=X^i\}$ and $\{b_j=Y^j\}$,
whence $A\otimes B=\F[X,Y]/(f(X),g(Y))$.
\begin{prop}\label{prop:tensor}
  $S_n$ decomposes as
  $$S_n = R_{2^n-1}\oplus\F = \bigotimes\F[X]\big/\widehat{f_{n,i}}
  = \F[X_1,\dots,X_{r_n}]\big/\Big(\widehat{f_{n,1}}(X_1),\dots,\widehat{f_{n,r_n}}(X_{r_n})\Big).$$
\end{prop}

\begin{cor}\label{cor:tensor}
  If $\Phi_n$ factors in $r_n>1$ factors, then $\Phi_{2^n-1}$ factors in
  at least $2^{r_n}-1>1$ factors; if $f$ is a factor of $\Phi_n$, then
  $\widehat f/X$ is a factor of $\Phi_{2^n-1}$.
\end{cor}
\begin{proof}
  The factors $f_{n,i}$ of $\Phi_n$ are irreducible, but the
  $\widehat{f_{n,i}}$ have at least two factors, one
  of them being $X$. According to the proposition,
  $$S_n=R_{2^n-1}\oplus\F=\bigotimes\F[X]/\widehat{f_{n,i}}=\bigotimes\left(\F[X]\big/(\widehat{f_{n,i}}/X)\oplus\F\right).$$
  If we distribute the $r_n$ direct sums over the tensor products, we
  obtain an expression of $S_n$ as a direct sum of $2^{r_n}$ algebras.
  Among these is $\F=\F\otimes\dots\otimes\F$; all the $2^{r_n}-1$
  others are summands of $R_{2^n-1}$. Among these others are the
  $\F\otimes\dots\otimes\F[X]/(\widehat{f_{n,i}}/X)\otimes\dots\otimes\F$.
\end{proof}

\begin{proof}[Proof of Proposition~\ref{prop:tensor}]
  By Lemma~\ref{lem:nrf}, $\Phi_n$ factors as claimed. By
  induction, it suffices to consider a factorization $\Phi_n=fg$, with
  $f$ and $g$ coprime, and to show that in that case
  $$S_n=\F[X]/\widehat{\Phi_n}\cong\F[X]/\widehat
  f\otimes\F[X]/\widehat g=\F[Y,Z]/(\widehat f(Y),\widehat g(Z)).$$
  As $f$ and $g$ are coprime, apply B\'ezout's theorem to decompose
  the identity $1=\alpha f+\beta g$, for $\alpha$ and $\beta$
  polynomials. Apply the ``hat'' operator:
  $$X = \widehat\alpha(\widehat f(X)) + \widehat\beta(\widehat g(X)).$$
  We may now define the two mutually inverse maps
  \begin{align*}
    \F[X]/(\widehat f\circ\widehat g) &\rightleftarrows \F[Y,Z]/(\widehat f(Y),\widehat g(Z))\\
    X &\rightarrow \widehat\beta(Y)+\widehat\alpha(Z)\\
    \widehat g(X) \leftarrow Y, & \phantom{\rightarrow{}}\widehat f(X) \leftarrow Z.
  \end{align*}
\end{proof}
Really, this proposition is a dual version of the Chinese Remainder
Theorem, and its proof draws largely on this fact: we constructed
natural injections $\F[X]/\widehat{f_{n,i}}\hookrightarrow S_n$ dual to
the natural projections $R_n\twoheadrightarrow\F[X]/f_{n,i}$.

The decomposition stated in Corollary~\ref{cor:tensor} need not be
complete, though, as the tensor product of fields need not be a
field:
\begin{prop}
  Let $f$ and $g$ be two irreducible polynomials, so that
  $A=\F[X]/f(X)$ and $B=\F[Y]/g(Y)$ are fields. Then $A\otimes
  B=\F[X,Y]/(f(X),g(Y))$ is a direct sum of $\gcd(\deg f,\deg g)$
  fields of degree $\lcm(\deg f,\deg g)$; in particular, $A\otimes B$
  is a field if and only if $\deg f$ and $\deg g$ are coprime.
\end{prop}
\begin{proof}
  $A\otimes B$ is semisimple and commutative whenever both $A$ and $B$
  are, so $A\otimes B$ is a direct sum of fields. Let $\chi_A$ and
  $\chi_B$ be the Frobenius automorphisms of $A$ and $B$: then the
  Frobenius automorphism of $A\otimes B$ is
  $\chi=\chi_A\otimes\chi_B$, so is of order exactly $\lcm(\deg f,\deg
  g)$, and all subfields of $A\otimes B$ are of degree at most
  $\lcm(\deg f,\deg g)$. On the other hand, $A\otimes B$ splits as a
  sum of fields each containing $A$ and $B$ (see~\cite[page
  54]{frohlich:ant}).
\end{proof}

We give an example of Corollary~\ref{cor:tensor} in the first
non-trivial case, $n=5$: then $\Phi_5 = X^5+X+1 = (X^2+X+1)(X^3+X^2+1)$, so
$$R_n = \F[X]/(X^2+X+1) \oplus \F[Y]/(Y^3+Y^2+1).$$
Let us note $f = X^2+X+1$, $g = Y^3+Y^2+1$, and for convenience
$F=\widehat f/X=X^3+X+1$ and $G=\widehat g/Y=Y^7+Y^3+1$. Then
\begin{align*}
  S_n &= \F[X]/\widehat f \otimes \F[Y]/\widehat g\\
      &= (\F\oplus\F[X]/F) \otimes (\F\oplus\F[Y]/G)\\
      &= \F \oplus \F[X]/F \oplus \F[Y]/G \oplus \F[X,Y]/(F,G)\\
      &= \F \oplus \F[X]/F \oplus \F[Y]/G \oplus \F[Z]/H,
\end{align*}
where $H(Z)=\Phi_{31}(Z)/F(Z)/G(Z)$.
This in turn factors
$$\Phi_{31}=\widehat f/X\cdot\widehat g/X\cdot \left(X\frac{\widehat f\circ \widehat g}{\widehat f\widehat g}\right).$$
Note that the factors are not necessarily irreducible.

Finally there is an interesting connection between the orbits under
multiplication by $X$ in $R_n$ and the factorization of $\Phi_{2^n-1}$.
\begin{prop}
  Let $\mathcal O_0=\{0\}$ and $\mathcal O_i$ for $i\in\{1,\dots,k\}$
  be the orbits in $R_n$ under the ``multiply-by-$X$'' $\F[X]$-action.
  Then $|\mathcal O_i|=\ell=(2^n-1)/k$ for all $i$, and $S_n$ splits
  as $\F\oplus A\oplus\dots\oplus A$ (with $k$ copies of $A$), where
  $A$ is an algebra of dimension $\ell$.
\end{prop}
\begin{proof}
  A reformulation of Proposition~\ref{prop:2^n-1} is that there is a
  group homomorphism $R_n\supset\langle
  X\rangle\to\operatorname{Aut}(S_n)$, mapping $X$ to $\chi$. Now by
  assumption $X$ is of order $\ell$ in $R_n$, so $\chi^\ell=1$ in
  $S_n$, and $S_n$ splits as a direct sum of fields of degree dividing
  $\ell$.

  A generator of $R_n^\times$ maps to an automorphism of
  $\widehat{R_n}$, whose $k$th power is $\chi$. It must act by
  permutation and automorphisms on a set of $k$ subalgebras of $S_n$,
  who are then all isomorphic; call them $A$.
\end{proof}
The first values of $n$ for which $R_n$ is a field are
$2,3,4,6,7,9,15,22,28,30,46$. In the first non-trivial example,
$\Phi_9$ is irreducible, but $\Phi_{511}$ is the product of seven
polynomials of degree $73$.

We are now naturally led to the following
\begin{conj}
  Let $\mathcal O_i$ for $i\in\{0,\dots,k\}$ be the orbits in
  $R_n$ under the ``multiply-by-$X$'' $\F[X]$-action. Then $S_n$
  splits as a direct sum of fields of degree $|\mathcal O_i|$
  for all $i\in\{0,\dots,k\}$.
\end{conj}
This conjecture generalizes and contains Conjecture~\ref{conj:2^n-1}.

\def\nop#1{}
\providecommand{\bysame}{\leavevmode\hbox to3em{\hrulefill}\thinspace}

\end{document}